%Apanasov's paper on bending deformations

\input amstex
\input psfig
\documentstyle{amsppt}
\magnification=1200
\NoBlackBoxes
\hcorrection{0.5truein}   %?
\TagsOnRight
\nologo
\define\Ga{\Gamma}
\define\ga{\gamma}
\define\Om{\Omega}
\define\La{\Lambda}
\define\p{\partial}
\define\a{\alpha}
\define\sa{\sigma}

\define\da{\delta}

\define\e{\epsilon}

\define\vp{\varphi}
\define\vP{\varPhi}

\define\sch{{\Cal H}}

\define\sct{{\Cal T}}
\define\ba{{\Bbb A}}

\define\bc{{\Bbb C}}

\define\bh{{\Bbb H}}

\define\br{{\Bbb R}}

\define\gtS{{\frak S}}
\define\gts{{\frak s}}

\def\ref#1{[#1]}
\define\ch{\bh_{\bc}^}
\define\rh{\bh_{\br}^}

\define\ra{\rightarrow}
\define\hra{\hookrightarrow}
\define\bs{\backslash}
\define\ov{\overline}
\define\col{\!:\!}
\define\Hom{\operatorname{Hom}}

\define\Vol{\operatorname{Vol}}

\define\is{\operatorname{Isom}}

\define\arccsch{\operatorname{arccsch}}
\define\csch{\operatorname{csch}}

\parindent=1pc                  %% AMSPPT default setting

\leftheadtext{Boris Apanasov}
\rightheadtext{Bending deformations}
\topmatter
\title Bending deformations of complex hyperbolic surfaces
\endtitle

\author Boris Apanasov\footnote"\dag"{ 
Supported in part by the NSF.
%Research at MSRI was supported in part by NSF grant DMS-9022140.
\hfill \hfill \hfill}
\endauthor
\address 
 Dept. of Math., Univ. of Oklahoma,
Norman, OK  73019  \hfill \hfill \hfill \endaddress
\address 
 Mathematical Sciences Research Institute, 
%1000 Centennial Dr.,
Berkeley, CA 94720-5070 \hfill \hfill \hfill \endaddress
\email apanasov\@uoknor.edu\hfill \hfill \hfill\endemail

%\affil (Preliminary version)\endaffil
\subjclass 57,55,53,51,30
\endsubjclass
\abstract\nofrills {ABSTRACT.}
We study deformations of complex hyperbolic surfaces which furnish the simplest
examples of : \,(i) negatively curved K\"ahler manifolds and 
(ii) negatively curved Riemannian manifolds not having {\it constant} curvature.
Although such complex surfaces may share the rigidity of
quaternionic/octionic hyperbolic manifolds, our main goal is to show
that they enjoy nevertheless the flexibility of low-dimensional real hyperbolic manifolds.
Namely we define a class of ``bending" deformations of a given (Stein) complex surface $M$ 
associated with its closed geodesics provided that $M$ is homotopy equivalent to a Riemann 
surface whose embedding in $M$ has a non-trivial totally real geodesic part.
Such bending deformations bend $M$ along its closed geodesics and are induced by equivariant 
quasiconformal homeomorphisms of the complex hyperbolic space and its Cauchy-Riemannian 
structure at infinity.
 \endabstract
\endtopmatter
\bigskip

\document
\head 1. Introduction \endhead

A complex hyperbolic manifold $M$ is locally modeled on the complex
hyperbolic space $\ch n$ which can be represented by the unit complex
ball $\Bbb B^n\subset \bc^n$ with the K\"ahler structure given by the
Bergman metric (with pinched negative curvature $K$, $-1\leq K \leq -
1/4$, compare \cite{M2, GP, G3, AX1}). In particular, in complex dimension two, 
due to Yau's uniformization theorem \cite{Ya}, every smooth complex
projective surface $M$ with positive canonical bundle satisfying the
topological condition $\chi(M)=3\cdot {\text {\it Signature}}(M)$ is a complex
hyperbolic manifold.
 Since $M=\ch n/G$, where the discrete
torsion free group $G\subset PU(n,1)$ is the fundamental group of $M$,
one can reduce the study of the Teichm\"uller space  $\Cal T(M)$ of isotopy 
classes of complex hyperbolic structures on $M$ to studying the variety
$\Cal T(G)$ of conjugacy classes of discrete faithful representations 
$\rho\col G\ra PU(n,1)$ (involving the space $\Cal D(M)$ of the developing maps, see
\cite{G2, FG}). Here $\Cal T(G)=\Cal R_0(G)/PU(n,1)$, and the variety 
$\Cal R_0(G)\subset \Hom (G, PU(n,1))$ consists of discrete faithful 
representations $\rho$ of the group $G$ whose co-volume, $\Vol (\ch n/G)$, 
may be infinite.

Due to the Mostow rigidity theorem \cite{M1}, hyperbolic structures of
finite volume and (real) dimension at least three are uniquely determined by
their topology, and one has no continuous deformations of them. Despite
that, real hyperbolic manifolds $M$ can be deformed as conformal manifolds,
or equivalently as higher-dimensional hyperbolic manifolds $M\times (0,1)$ 
of infinite volume. First such deformations were given by the author \cite{A2} and,
after Thurston's ``Mickey Mouse" example \cite{Th}, were called bendings of 
$M$ along its totally geodesic hypersurfaces, see also \cite{A1, A3-A5, JM, Ko}.
Furthermore all these deformations are quasiconformally equivalent showing 
a rich supply of quasiconformal $G$-equivariant homeomorphisms in the closure
of the real hyperbolic space $\rh n$. In particular the limit set 
$\La(G)\subset\p\rh {n+1}$ deforms continuously from a round sphere 
$\p \rh n=S^{n-1}\subset S^n=\rh {n+1}$ into a nondifferentiably embedded topological 
$(n-1)$-sphere quasiconformally equivalent to $S^{n-1}$.

Contrasting to the above flexibility, ``non-real" hyperbolic manifolds
seem much more rigid. In particular, due to Pansu \cite{P}, quasiconformal maps in the
sphere at infinity of quaternionic/octionic hyperbolic spaces are
necessarily automorphisms, and thus there cannot be interesting
quasiconformal deformations of corresponding structures. Secondly, due 
to Corlette's rigidity theorem \cite{Co}, such manifolds are even super-rigid
-- analogously to Margulis super-rigidity in higher rank \cite{Ma}.
Furthermore, complex hyperbolic manifolds share the above rigidity of  
quaternionic/octionic hyperbolic manifolds. Namely, due to the Goldman's
local rigidity theorem in dimension $n=2$ \cite{G1}, every nearby discrete representation 
$\rho\col G\to PU(2,1)$ of a cocompact lattice $G\subset PU(1,1)$ stabilizes a complex
geodesic in the complex hyperbolic space $\ch 2$ (which is also true for small deformations
of cocompact lattices $G\subset PU(n-1,1)$ in higher dimensions $n\geq 3$
\cite{GM}), and thus the limit set $\La(\rho G)\subset \p\ch n$
is always a round sphere $S^{2n-3}$. In higher dimensions $n\geq 3$, this local rigidity 
of complex hyperbolic $n$-manifolds $M$ homotopy equivalent to their closed complex 
totally geodesic hypersurfaces is even global due to a recent Yue's theorem \cite{Yu}.
  
Our goal here is to show that, in contrast to rigidity of complex hyperbolic 
n-manifolds $M$ from the above class, complex hyperbolic (Stein) manifolds $M$ 
homotopy equivalent to their closed totally 
{\it real} geodesic surfaces are not rigid. Namely, in complex dimension 2, we provide a canonical 
construction of continuous non-trivial quasi-Fuchsian deformations of complex surfaces
fibered over closed Riemannian surfaces. This is the first such
deformations of fibrations with compact base (for non-compact base, see
a different Goldman-Parker' deformation \cite{GP} of ideal triangle groups
$G\subset PO(2,1)$).
Our construction is inspired by the well know bending deformations 
of real hyperbolic (conformal) manifolds along totally geodesic hypersurfaces
and by a M.Carneiro--N.Gusevskii' construction of a 
discrete representation in $PU(2,1)$ of a surface group, 
which has been kindly showed us by N.Gusevskii \cite {Gu}. 
In the case of complex hyperbolic (and Cauchy-Riemannian) structures, it works 
however in a different way than in the real one.
Namely our complex bending deformations involve
simultaneous bending of the base of the fibration of the complex surface $M$ as well 
as bendings of each of its totally geodesic fibers (see Remark 4.10). 
Such bending deformations of complex surfaces are associated
to their real simple closed geodesics (of real codimension 3), but have
nothing common with cone deformations of real hyperbolic 3-manifolds along
closed geodesics (see \cite{A4, A5}). 

We also remark that the condition that the group $G\subset PU(n,1)$ preserves a complex
totally geodesic hyperspace in $\ch n$ is essential for local rigidity
of deformations only for co-compact lattices $G\subset PU(n-1,1)$. As
our subsequent work \cite{ACG} shows, there are non-trivial quasi-Fuchsian
deformations of such a {\it co-finite} lattice  $G\subset PU(1,1)$ in
$PU(2,1)$ induced by a continuous family of $G$-equivariant homeomorphisms 
in $\ov{\ch 2}$, which cannot however be quasiconformal due to some obstruction.

Furthermore, there are well known complications (cf. \cite{KR, P, V})
in constructing equivariant homeomorphisms in the complex hyperbolic space and 
in Cauchy-Riemannian geometry, which are due to necessary conditions for
such maps to preserve the K\"ahler and contact structures (correspondingly 
in the complex hyperbolic space and at its infinity, the one-point
compactification of the Heisenberg group$\sch_n$). Despite that, as it follows from our 
construction, the complex bending deformations are induced by equivariant 
homeomorphisms which are in addition quasiconformal with respect to the corresponding 
metrics. One of our main results may be formulated as follows.
 
\proclaim{Theorem 4.7} Let $G\subset PO(2,1)\subset PU(2,1)$ be a given
lattice uniformizing a Riemann 2-surface $S_p$ of genus 
$p\geq 2$. Then, for any simple closed geodesic \break
$\a\subset S_p=H^2_{\br}/G$ and a sufficiently small $\eta_0>0$, 
there is a bending deformation \break
$\Cal B_\a\col (-\eta_0,\,\eta_0)\ra \Cal R_0(G)$ of the group $G$ along $\a$,
$\Cal B_\a(\eta)=\rho_{\eta}=F^*_{\eta}$, 
induced by $G$-equivariant quasiconformal homeomorphisms 
$F_{\eta}: \ov{\ch 2} \ra \ov{\ch 2}$. 
\endproclaim

We notice that such complex bending deformations depend on many independent
parameters, as it is shown by application of our construction and 
\'Elie Cartan \cite {Ca} angular invariant in Cauchy-Riemannian
geometry:

\proclaim{Corollary 4.8} Let $S_p=\rh 2/G$ be a closed totally real
geodesic surface of genus $p>1$ in a given complex hyperbolic surface
$M=\ch 2/G$, $G\subset PO(2,1)\subset PU(2,1)$. Then there is an
embedding $\pi\circ\Cal B\col B^{2p-2}\hra \Cal T(M)$ of a real $(2p-2)$-ball
into the Teichm\"uller space of $M$,
defined by bending deformations along disjoint closed geodesics in $M$  and
the projection $\pi\col \Cal D(M)\ra \Cal T(M)=\Cal D(M)/PU(2,1)$.
\endproclaim
 
In our subsequent work \cite{AG}, we apply the constructed bending deformations
to answer a well known question about cusp groups on the boundary of
the Teichm\"uller space of $\sct (M)$ of a Stein complex surface $M$
fibering over a compact Riemann surface of genus $p>1$:
  
\proclaim {Theorem } Let $G\subset PO(2,1)\subset PU(2,1)$ be a
uniform lattice isomorphic to the fundamental group of a closed surface
$S_p$ of genus $p\geq 2$. Then, for any simple closed geodesic
$\a\subset S_p=H^2_{\br}/G$, there is a continuous deformation $\rho_t=f^*_t$ 
induced by $G$-equivariant quasiconformal homeomorphisms 
$f_t: \ov{\ch 2} \to \ov{\ch 2}$ whose limit representation $\rho_{\infty}$ corresponds
to a boundary cusp point of the Teichm\"uller space $\Cal T(G)$, that is
the boundary group $\rho_{\infty}(G)$ has an accidental parabolic element 
$\rho_{\infty}(g_{\a})$ where $g_{\a}\in G$ represents the geodesic $\a\subset S_p$.
\endproclaim

We note that, due to our construction of such continuous quasiconformal deformations
in \cite {AG},
they are independent if the corresponding geodesics $\a_i\subset S_p$
are disjoint. It implies the existence of a boundary group in $\p \Cal T(G)$ with
``maximal" number of non-conjugate accidental parabolic subgroups:

\proclaim {Corollary } Let $G\subset PO(2,1)\subset PU(2,1)$ be a
uniform lattice isomorphic to the fundamental group of a closed surface
$S_p$ of genus $p\geq 2$. Then there is a continuous deformation 
$R\col \br^{2p-2}\to \Cal T(G)$ whose boundary group
$G_{\infty}=R(\infty)(G)$ has $2p-2$ non-conjugate accidental parabolic
subgroups.
\endproclaim

Finally, we mention another aspect of the intrigue problem on geometrical finiteness
of complex hyperbolic surfaces (see \cite {AX1, AX2}) for which it may be possible to 
apply our complex bending deformations:

\proclaim {Problem } Construct a geometrically infinite (finitely
generated) discrete group $G\subset PU(2,1)$ whose limit set is the whole sphere
at infinity, $\La(G)=\p \ch 2=\ov {\sch}$, and which is the limit of
convex cocompact groups $G_i\subset PU(2,1)$ from the Teichm\"uller space 
$\Cal T(\Ga)$ of a convex cocompact group $\Ga\subset PU(2,1)$. Is that
possible for a Schottky group $\Ga$?
\endproclaim
\noindent{\it Acknowledgements.} The author would like to thank Nikolay
Gusevskii for helpful conversations. A part of paper was written during
the author's stay at MSRI. Research at MSRI was 
supported in part by NSF grant DMS-9022140.

\head 2. Complex hyperbolic surfaces homotopy equivalent to a Riemann surface
 \endhead

Let $M$ be a complex hyperbolic surface (a complete K\"ahler 2-manifold of 
constant holomorphic sectional curvature -1) homotopy equivalent to a Riemann
surface $S_g$ of genus $g>1$, that is $M=\ch 2/G$ is a quotient of the
complex hyperbolic 2-space $\ch 2$ by a discrete torsion free isometry
group $G\subset PU(2,1)$ isomorphic to the fundamental group
$\pi_1(S_g)$. The boundary at infinity of $M$ is a 3-dimensional manifold
with the induced spherical Cauchy-Riemannian structure modeled on the one-point 
compactification $\ov\sch=\sch\cup\{\infty\}\approx S^3$ of the Heisenberg 
group $\sch=\bc\times\br$.

We are interested in Teichm\"uller spaces of such complex
surfaces and Cauchy-Riemannian 3-manifolds or equivalently, 
in varieties of conjugacy classes of discrete
faithful representations $\rho\col G\to PU(2,1)$, especially in curves
in these spaces corresponding to continuous families of quasiconformal
deformations of such structures or discrete groups.

As the simplest examples of such complex surfaces, one can take the
quotients $M_i=\ch 2/G_i$ corresponding to two embeddings of $\pi_1(S_g)$
as lattices acting on totally geodesic planes in $\ch 2$, 
either on complex geodesics (represented by $\ch 1\subset \ch 2$) 
or on totally real geodesic 2-planes (represented by $\rh 2\subset \ch
2$). In these cases the discrete surface group is either 
$G_1\subset PU(1,1)\subset PU(2,1)$
or $G_2\subset PO(2,1)\subset PU(2,1)$, correspondingly. The limit
sets $\La(G_1)$ and $\La(G_2)$ are homeomorphic to the circle $S^1$. Moreover,
the actions of the groups $G_1$ and $G_2$ on such circles are
(equivariantly) homeomorphically conjugate due to the following Apanasov 
\cite{A6} isomorphism theorem:

\proclaim{Theorem 2.1} Let $\varphi : G\ra H$ be a type preserving isomorphism
 of two non-ele\-men\-ta\-ry geometrically finite discrete subgroups 
$G,H\subset PU(n,1)$. Then
 there exists a unique equivariant ho\-meo\-mor\-phism 
$f_{\varphi}\col \La(G)\ra \La(H)$ of
 their limit sets that induces the isomorphism $\varphi$.
 \endproclaim
 
However, in the case of the above surface groups $G_1$ and $G_2$, the equivariant
homeomorphism $f_{\varphi}\col \La(G_1)\ra \La(G_2)$ cannot be homeomorphically 
extended to the whole sphere $\p \ch 2 \approx S^3$.
The obstruction to that is due to the fact that the quotient
complex surfaces $\ch 2/G_1$ and $\ch 2/G_2$ are not homeomorphic. Namely, these
complex surfaces are disc bundles over $S_g$ and have different Toledo
invariants: $\tau(\ch 2/G)=2g-2$ and $\tau(\ch 2/H)=0$, see \cite {To}.

The complex structures of the complex surfaces $M_1$ and $M_2$ are quite different, too.
While the manifolds $M_1$ (corresponding to $G_1\subset PU(1,1)$) have natural
embeddings of the Riemann surface $S_g$ as holomorphic totally geodesic submanifolds 
and hence cannot be Stein manifolds if $S_p\subset M_1$ is compact, the manifolds 
$M_2$  (corresponding to $G_2\subset PO(2,1)$) are Stein manifolds, see Burns-Shnider \cite{BS}. 
Moreover due to Goldman \cite{G1}, if the surface $S_p\subset M_1$ is closed 
(i.e. the lattices $G_1\subset PU(1,1)$ are co-compact), 
the manifolds $M_1$ are locally rigid in the sense that 
every nearby representation $G_1\rightarrow PU(2,1)$
stabilizes a complex geodesic in $\ch 2$ and is conjugate to a
representation $G_1\rightarrow PU(1,1) \subset PU(2,1)$. In other words,
there are no non-trivial ``quasi-Fuchsian" deformations of $G_1$ and $M_1$.

Our goal here is to show that, in contrast to such rigidity of manifolds $M_1$ 
obtained from uniform lattices in $PU(1,1)$, the Stein manifolds $M_2$
obtained from discrete subgroups in $PO(2,1)$ are non-rigid. Furthermore,
for a given simple geodesic $\a\subset M_2$ (in the totally real geodesic surface
$S_g\subset M_2$), there is a continuous non-trivial path in the Teichm\"uller
space $\sct (M_2)$ (or $\sct (G_2)$) represented by a continuous family of 
$G_2$-equivariant quasiconformal homeomorphisms $f_{\a,t}: \ov\sch\to\ov\sch $. 
We call such a deformation $\Cal B_{\a}$, $\Cal B_{\a}(t)=f_{\a,t}^* $,
 a {\it bending
deformation} of Cauchy-Riemannian 3-manifold  $\p_{\infty} M_2$, the group $G_2$ 
and the complex 
hyperbolic surface $M_2$ along a given closed geodesic $\a\subset M_2$.

\head 3. Colar lemma and Dirichlet polyhedra\endhead 

For a given discrete group $\Ga\subset PO(2,1)\subset PU(2,1)$, let $S=\rh 2/\Ga$ be a 
real hyperbolic 2-orbifold embedded in $M=\ch 2/\Ga$ as a totally real
geodesic suborbifold. While $M$ has constant holomorphic sectional
curvature  -1, its sectional curvature varies in [-1, -1/4], and in particular 
the real 2-orbifold $S$
has constant sectional curvature $K=-1/4$. 

For a simple closed geodesic $\a\subset S\subset M$, let $\ga_{\a}\in \pi^{orb}_1(S)\cong \Ga$
represent it. Then the fundamental group
$ \pi^{orb}_1(S)\cong \Ga$ can be decomposed into either a free amalgamated
product or HNN-extension of its subgroups. Namely, for
$\Ga_0=<\ga_{\a}>$, we correspondingly have for separating and non-separating geodesics $\a$: 

$$ \Ga=\Ga_1 *_{\Ga_0}\Ga_2 \quad \text{and} \quad \Ga=\Ga_1 *_{\Ga_0}\,.$$

Here we would like to focus on an analogue of the well-known collar
lemma which (for Riemann surfaces) is originally due to L.Keen \cite {Ke} and 
its sharp form due to P.Buser \cite{Bu}. It claims that a simple closed
geodesic $\a$ of length $\ell$ on a Riemann surface $S$ of
constant curvature -1 has a collar neighborhood $U_{\da}(\a)\subset S$
of radius $\da=\ln (\coth \ell/4)$, that is no disjoint closed geodesic $\beta\subset S$, 
$\a\cap\beta=\emptyset$, intersects the neighborhood $U_{\da}(\a)$.

To reformulate the collar lemma for isometric action of a discrete group
$\Ga\subset PO(2,1)\subset PU(2,1)$
on the real hyperbolic plane $\rh 2\subset \ch 2$, we denote $A\subset \rh 2$ the axis of
a hyperbolic element $\ga_{\a}\in\Ga$ (the translation along $A$ of the
length $\ell $). Namely, multiplying lengths by $\sqrt {-K}$, we have due to 
Beardon \cite {Be, 11.6.10} that, for any such discrete
isometric $\Ga$-action on the real hyperbolic plane of constant negative
curvature $K<0$,

$$\sinh \left(\frac{\ell \sqrt {-K}}{2}\right)\cdot\sinh
 \left(\frac{d(A,g(A)) \sqrt {-K}}{2}\right) \geq
\frac{1}{2}\,,$$
for any $\ga\in\Ga\bs\Ga_0$, $\,\,\Ga_0=<\ga_{\a}>$.

This shows that, for any element $\ga\in\Ga\bs\Ga_0$ of the surface group 
$\Ga\subset PO(2,1)\subset PU(2,1)$, $\Ga\cong \pi_1(S)$, which does not preserve 
the axis $A$ of the cyclic hyperbolic group $\Ga_0=<\ga_{\a}>$,

$$\sinh \left(\frac{\ell}{4}\right)\cdot\sinh
 \left(\frac{d(A,g(A))}{4}\right) \geq \frac{1}{2}\,.\tag3.1$$

\proclaim{Lemma 3.2}
 Let $\ga_\a$ be a hyperbolic element of a non-elementary discrete group 
$\Ga\subset PO(2,1)\subset PU(2,1)$ with translation length $\ell$ along its axis 
$A\subset \rh 2$, and $\da>0$ a radius of a tubular neighborhood $U_{\da}(A)$ precisely 
invariant with respect to its stabilizer $\Ga_0\subset \Ga$:

$$ \Ga_0(U_\da(A))= U_\da(A) \quad \text{and} \quad \ga((U_\da(A))\cap U_\da(A)
=\emptyset, \quad \ga\in \Ga\bs\Ga_0\,.$$

Then the Dirichlet polyhedron
$D_z(\Ga)$ of the group $\Ga$ centered at a point $z\in A$ has two sides $\sa, \sa'$, 
$\ga_\a(\sa)=\sa'$, provided that $\da>\ell/4$.
\endproclaim

\demo{Proof} Since the central point $z$ lies on the $\Ga$-invariant
totally real geodesic plane $\rh 2\subset \ch 2$, the Dirichlet polyhedron
$D_z(\Ga)$ of the group $\Ga$ has the same combinatorics as its
restriction to $\rh 2$ (the Dirichlet polygon there), see \cite {G3,
VIII.4}.  On the plane $\rh 2$, the worst case is the case of an element
$\ga\in \Ga\bs\Ga_0$ that maps the point $z$ to a point $\ga(z)$ which
lies in the intersection of one of Dirichlet bisectors 
$\frak S(z,\ga_\a(z))$ and $\frak S(z,\ga_\a^{-1}(z))$,

$$ \frak S(z_1, z_2)=\left\{ w\in \ch 2\,\col\, d(z_1,w)=d(z_2,w)\right\} \,,$$
and the boundary $\p U_{2\da}(A)$ of the 
$2\da$-neighborhood of the axis $A$. Considering the right geodesic triangle
with vertices $z, \ga(z)$ and $A\cap\frak S(z,\ga_\a(z))$, we clearly see that
the bisector $\frak S(z,\ga(z))$ passes the point  $A\cap\frak S(z,\ga_\a(z))$
if and only if $\ell=4\da$. This completes the proof. \qed
\enddemo

\proclaim{Corollary 3.3}
Let $\ga_\a$ be a hyperbolic element of a non-elementary discrete group 
$\Ga\subset PO(2,1)\subset PU(2,1)$ with translation length $\ell$ 
along its axis $A\subset \rh 2$. Then any tubular neighborhood 
$U_{\da}(A)$ of the axis $A$ of radius $\da>0$  is precisely 
invariant with respect to its stabilizer $\Ga_0\subset \Ga$ if

$$ \sinh \left(\frac{\ell}{4}\right)\cdot\sinh \left(\frac{\da}{2}\right)\leq 
\frac{1}{2}\,.\tag 3.4$$

Furthermore, for sufficiently small $\ell$, $\ell<4\da$, the Dirichlet 
polyhedron $D_z(\Ga)$ of the group $\Ga$ centered at a point $z\in A$ 
has two sides $\sa$ and $\sa'$ intersecting the axis $A$ and such that 
$\ga_\a(\sa)=\sa'$.
\endproclaim

\remark{Remark 3.5} It is easy to see that both conditions of 
Corollary hold for $\ell>0$ and $\da>0$ such that 

$$ \gather 
\sinh \frac{\ell}{2} < \frac{1}{2} \csch \left(\frac{\ln 63}{4}\right)
\approx 0.406114\\
\ln 63\leq \da < \arccsch \left( 2\sinh \frac{\ell}{4}\right)
\,. 
\endgather $$
The first condition in particular holds for $\ell<0.791411$ (see also \cite{GKL}).
\endremark

\demo{Proof} The claim about the sides of the Dirichlet polyhedron $D_z(\Ga)$
directly follows from the Lemma 3.2. In particular, it is easy to verify that
the values of the length $\ell$ and the
radius $\da$ in Remark 3.5 satisfy the inequality $\ell < 4\da$.

 To show that $U_{\da}(A)$ is precisely $\Ga_0$-invariant in $\ch 2$ it 
is enough to prove this fact for the restriction of $\Ga$-action to the
plane $\rh 2$ (see also \cite{GKL}). Here, due to (3.1) and the hypothesis (3.4), we have:

$$\sinh \frac{\da}{2}\leq\frac{1}{2}\csch \left(\frac{\ell}{4}\right)\leq 
\sinh \left(\frac{d(A,g(A))}{4}\right)\,.$$

This shows that 

$$\sinh \frac{\da}{2} < \sinh \frac{d(A,\ga(A))}{4} \quad \text {for
any} \quad \ga\in \Ga\bs \Ga_0\,,$$
which implies that $d(A,\ga(A))>2\da$, and hence
$ \ga((U_\da(A))\cap U_\da(A)
=\emptyset$ for any $\ga\in \Ga\bs\Ga_0$. \qed
\enddemo 

\head 4. Bendings of complex surfaces along simple closed geodesics
      \endhead
 
Now we start with a totally real geodesic surface $S=\rh 2/\Ga$ in the complex surface
$M=\ch 2/\Ga$, where $\Ga\subset PO(2,1)\subset PU(2,1)$ is a
given discrete group, and fix a simple geodesic $\a$ on $S$. According to the previous
sections 2 and 3, we may assume that the loop $\a$ is covered by a geodesic 
$A\subset \rh 2\subset \ch 2$ whose ends at infinity are $\infty$ and the origin
of the Heisenberg group $\sch=\bc \times \br$, $\ov \sch =\p \ch 2$.
Furthermore, deforming the surface $S$ and  its holonomy group
$\Ga\subset PO(2,1)$ (bending along the geodesic $\a$), we can assume
that the hyperbolic length of $\a$ is sufficiently small and the radius of its 
tubular neighborhood is big enough - as in Corollary 3.3 and Remark 3.5.
 
Then the group $G$ and its subgroups $G_0, G_1, G_2$ in the free
amalgamated (or HNN-extension) decomposition of $G$ have Dirichlet
polyhedra $D_z(G_i)\subset \ch 2$, $i=0,1,2$, centered at a point $z\in
A=(0,\infty)$, whose intersections with the hyperbolic 2-plane $\rh 2$
have the shapes indicated in Figures 1-4.

\midinsert
$$\hss\psfig{file=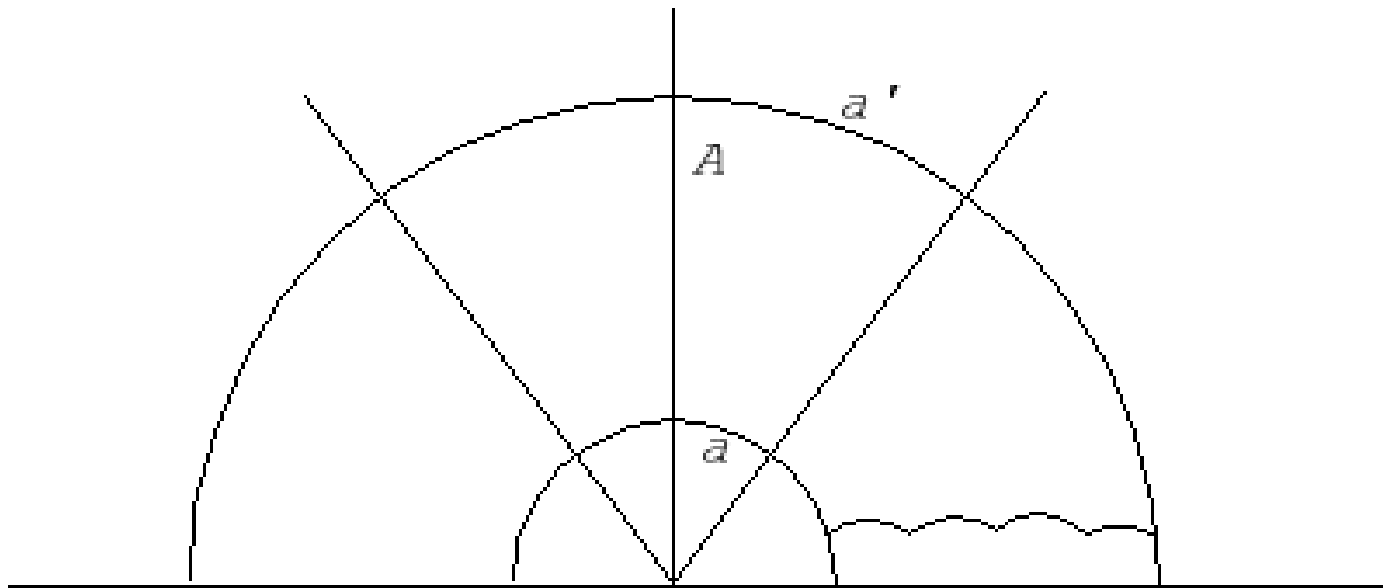}\hss$$
\botcaption{Figure 1} 
$G_1\subset G=G_1 *_{G_0}G_2$.
\endcaption
\endinsert

\midinsert
$$\hss\psfig{file=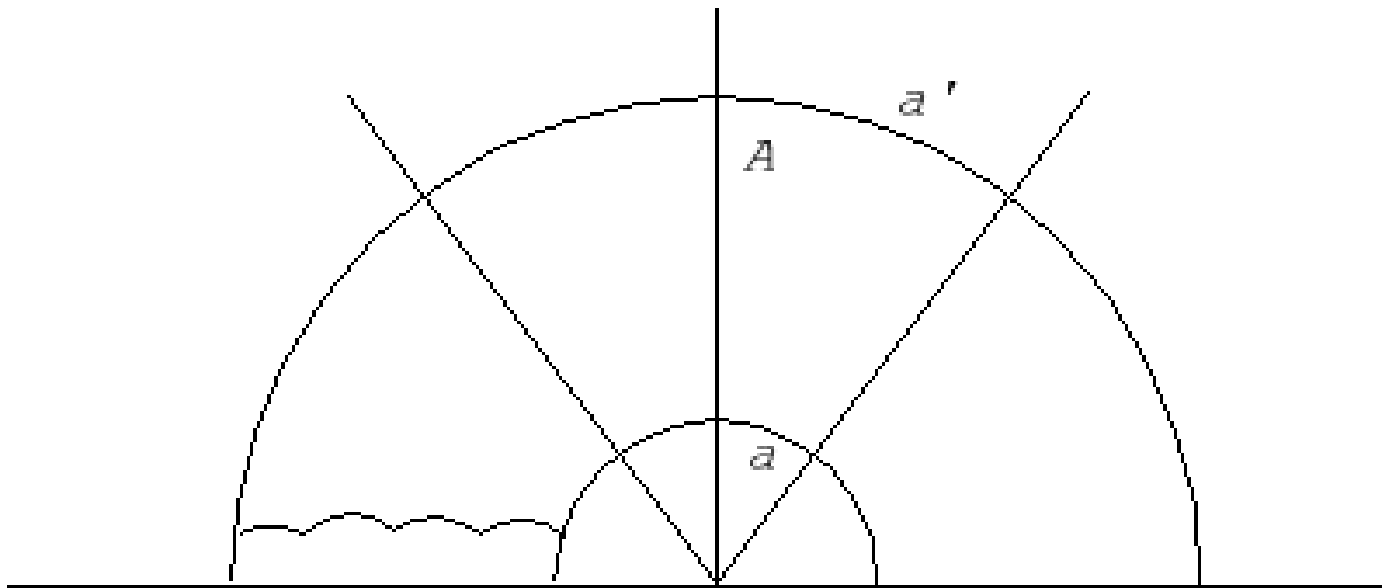}\hss$$
\botcaption{Figure 2} 
$G_2\subset G=G_1 *_{G_0}G_2$.
\endcaption
\endinsert

\midinsert
$$\hss\psfig{file=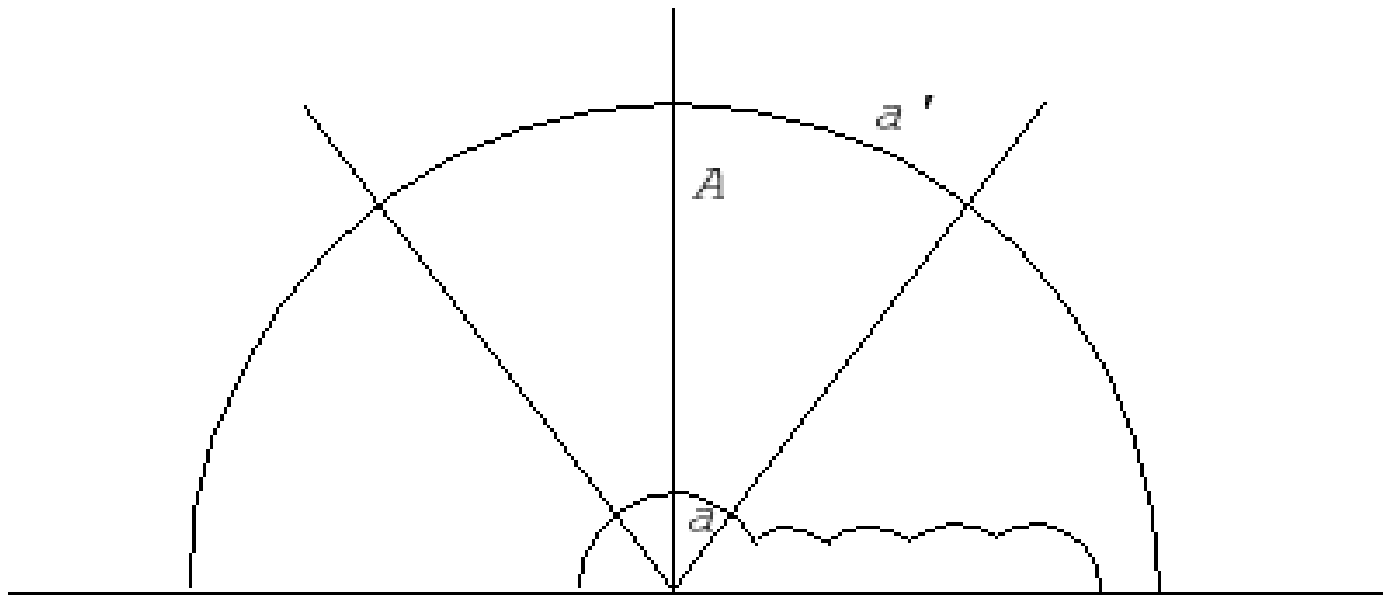}\hss$$
\botcaption{Figure 3} 
$G_1\subset G=G_1 *_{G_0}$.
\endcaption
\endinsert

\midinsert
$$\hss\psfig{file=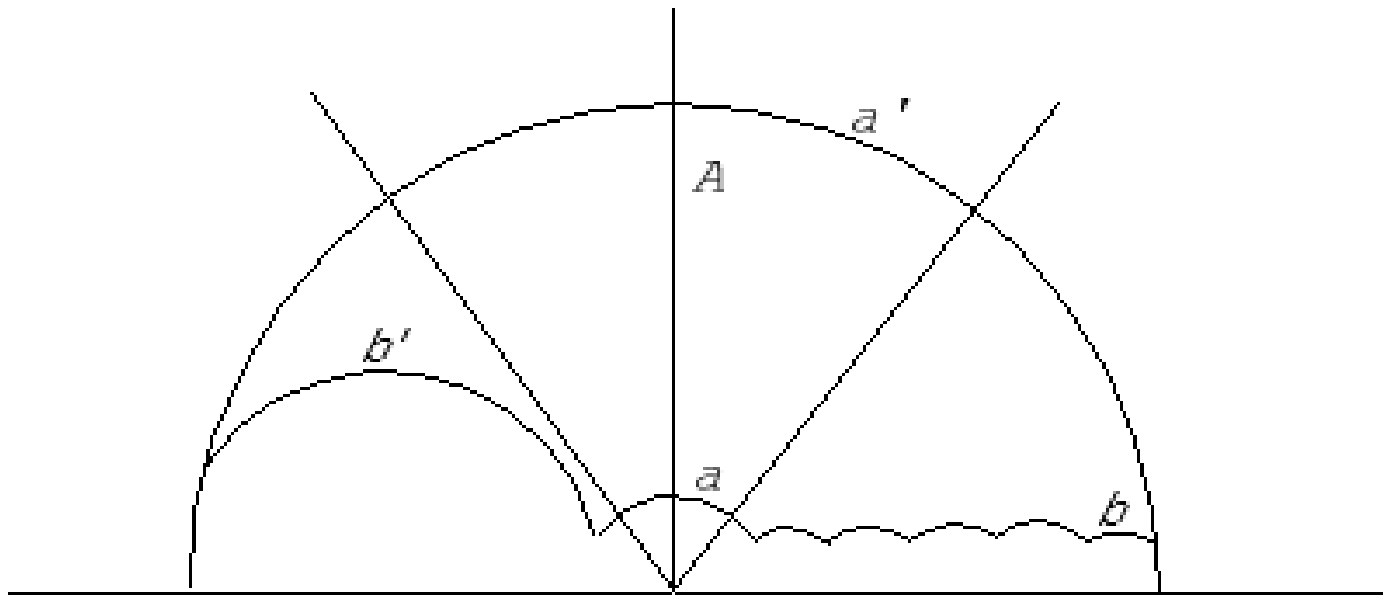}\hss$$
\botcaption{Figure 4} 
$G=G_1 *_{G_0}$.
\endcaption
\endinsert

In particular we have that, except two bisectors $\frak S$ and $\frak S'$ that are 
equivalent under 
the hyperbolic translation $g_{\a}$ (which generates the stabilizer $G_0$ of the 
axis $A$),
all other bisectors bounding those Dirichlet polyhedra lie in sufficiently small ``cone
neighborhoods" $C_+$ and $C_-$ of the arcs (infinite rays) $\br_+$ and $\br_-$ of the real circle
$\br\times \{0\}\subset \bc\times \br=\sch$. 

Actually, we may assume that the
Heisenberg spheres at infinity of the bisectors  $\frak S$ and $\frak S'$ 
have radii 1 and $r_0>1$, correspondingly. Then, for a
sufficiently  small $\e$, $0<\e<<r_0-1$, the cone neighborhoods 
$C_+, C_-\subset \ov {\ch 2}\bs\{\infty\}= \bc\times\br\times [0,+\infty)$
are correspondingly the cones of the $\e$-neighborhoods of the points
$(1,0,0), (-1,0,0) \in \bc\times\br\times [0,+\infty)$ with respect to the 
Cygan metric $\rho_c$ in $\ov {\ch 2}\bs\{\infty\}$. 
Here the Cygan metric $\rho_c$ is induced by the following norm (see \cite {Cy, Pr})
in the half-space model $\bc\times\br\times [0,\infty)$ of 
$\ov{\ch 2}\bs\{\infty\}$:

$$||(\xi,v,u)||_c=|\,||\xi||^2 + u - iv|^{1/2}\,.\tag4.1$$

 Clearly, we may
consider the length $\ell$ of the geodesic $\a$ so small that closures of all
equidistant halfspaces in $\ov {\ch 2}\bs\{\infty\}$ bounded by those bisectors and
disjoint from the Dirichlet polyhedron $D_z(G)$ do not intersect the co-vertical 
bisector whose infinity is $i\br\times\br\subset \bc\times \br$. It
follows from the fact (see \cite{G3, Thm VII.4.0.3}) that equidistant 
half-spaces  $\frak S_1$ and $\frak S_2$ in $\ch 2$ are disjoint if and
only if the intersection half-planes $\frak S_1\cap \rh 2$ and 
$\frak S_2\cap \rh 2$ are disjoint, see Figures 1-4.

Now we are ready to define a quasiconformal bending deformation of the group $G$
along the geodesic $A$, which defines a bending deformation of the complex
surface $M=\ch 2/G$ along the given closed geodesic $\a\subset S\subset
M$.

We specify numbers $\eta$ and $\zeta$ such that
$0<\zeta<\pi/2$, $0\leq\eta<\pi-2\zeta$ and the intersection $C_+\cap (\bc\times \{0\})$
is contained in the angle $\{z\in \bc\col |\arg z|\leq\zeta\}$. Then we
define a bending
homeomorphism $\phi=\phi_{\eta,\zeta}\col\bc\ra\bc$, see Fig. 5, which bends
the real axis $\br\subset\bc$ at the origin by the angle $\eta$:

$$\phi_{\eta,\zeta}(z)=\cases
z & \text{if}\ \, |\arg z|\geq\pi-\zeta\\
z\cdot\exp(i\eta) &  \text{if}\ \, |\arg z|\leq\zeta\\
z\cdot\exp(i\eta(1-(\arg z-\zeta)/(\pi-2\zeta))) & \text{if}\ \, \zeta
<\arg z<\pi-\zeta\\
z\cdot\exp(i\eta(1+(\arg z+\zeta)/(\pi-2\zeta))) & \text{if}\ \, \zeta
-\pi<\arg z<-\zeta\,.\endcases\tag4.2$$

\midinsert
$$\hss\psfig{file=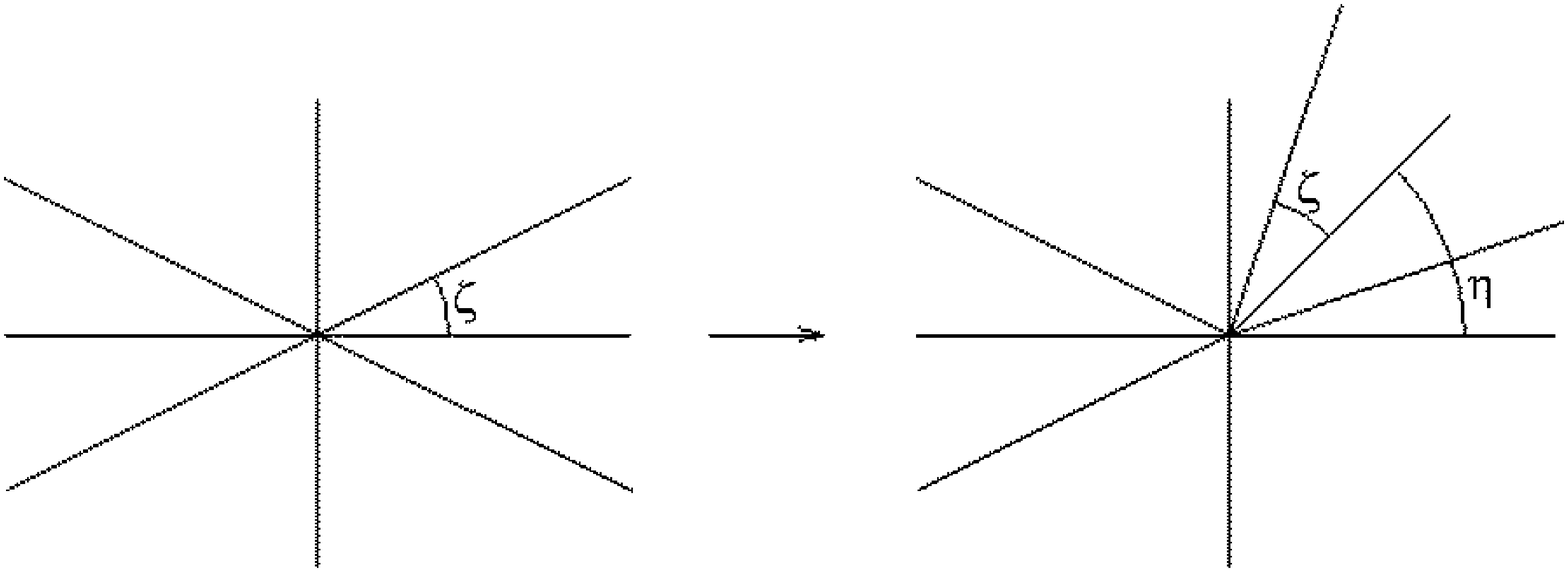,width=400pt}\hss$$
\captionwidth{7pc}
\botcaption{Figure 5} % usual picture for bending
\endcaption
\endinsert

For negative $\eta$, $2\zeta-\pi<\eta<0$, we set
$\phi_{\eta,\zeta}(z)=\ov{\phi_{-\eta,\zeta}(\ov z)}$.  Clearly,
$\phi_{\eta,\zeta}$ is quasiconformal with respect to the Cygan norm (4.1)
and is an isometry in the $\zeta$-cone
neighborhood of the real axis  $\br$ because its linear distortion is given by

$$K(\phi_{\eta,\zeta},z)=\cases
1\, & \text{if}\ \, |\arg z|\geq\pi-\zeta\\
1\, & \text{if}\ \, |\arg z|\leq\zeta\\
(\pi-2\zeta)/(\pi-2\zeta-\eta)\, & \text{if}\ \, \zeta<\arg z<\pi-\zeta\\
(\pi-2\zeta+\eta)/(\pi-2\zeta)\, & \text{if}\ \, \zeta-\pi<\arg z<-
\zeta\,.\endcases\tag4.3$$

Foliating the punctured Heisenberg group $\sch\bs\{0\}$ by Heisenberg
spheres $S(0,r)$ of radii $r>0$, we can extend this bending
homeomorphism to the whole sphere $S^3=\ov{\sch}$ at infinity. 

Namely,
let $W_+, W_-\subset\sch$ be ``dihedral angles" with the common vertical
axis $\{0\}\times \br$, which are foliated by arcs of real circles
connecting points $(0,v)$ and $(0,-v)$ on the vertical axis and
intersecting the the $\zeta$-cone neighborhoods of infinite rays
$\br_+, \br_- \subset \bc$, correspondingly. For each Heisenberg sphere
$S(0,r)$, there is a unique projective (``conformal" with respect to the Cygan metric) 
transformation $h_r\in PU(2,1)$ that maps $S(0,r)$ to the horizontal sphere
$(\bc\times \{0\})\cup \{\infty\}$ such that:

\roster
\item  It maps the intersection points 
$S(0,r)\cap (\{0\}\times\br)$ to the origin and infinity;

\item The chain $S(0,r)\cap (\bc\times \{0\})$ is pointwise fixed by $h_r$;

\item The intersections $S(0,r)\cap W_+$ and  $S(0,r)\cap W_+$ 
are mapped onto the $\zeta$-cone neighborhoods of 
$\br_+, \br_- \subset \bc$, correspondingly.
\endroster

Now we define an elementary bending homeomorphism 
$\vp=\vp_{\eta,\zeta}\col\sch\ra\sch$
as the homeomorphism whose restrictions to any Heisenberg sphere
$S(0,r)$ are conjugations of the bending $\phi_{\eta,\zeta}$ in the
horizontal plane $\bc\times\{0\}$ by the projective transformations $h_r$,

$$ \vp_{\eta,\zeta}\col\sch\ra\sch\,,\quad
\vp_{\eta,\zeta}\big |_{S(0,r)}= h_r^{-1}\phi_{\eta,\zeta}h_r\,.\tag4.4$$

Since the transformations $h_r$ map chains in $S(0,r)$ (horizontal
circles) to chains in the horizontal plane $\bc\times\{0\}$ (circles
centered at the origin), it directly follows from our construction 
that the defined homeomorphism
$\vp_{\eta,\zeta}$ preserves the horizontal plane $\bc\times\{0\}$, and
its restriction to this plane coincides with the bending $\phi_{\eta,\zeta}$. 
In other words, the elementary bending $\vp_{\eta,\zeta}$ bends the
plane $\br\times\br\subset \sch$ along the vertical axis
$\{0\}\times\br$, by a given angle $\eta$. Defining $\vp(\infty)=\infty$
and $\vp(0)=0$, 
we naturally extend this homeomorphism to the whole sphere $\ov{\sch}$. 

It follows from 
(4.3) that $\vp_{\eta,\zeta}$ is a $G_0$-equivariant quasiconformal homeomorphism 
in $\ov{\sch}$. Moreover, its restrictions to the ``dihedral angles" $W_-$ 
and $W_+$ are correspondingly the identity and the unitary rotation
$U_{\eta}\in PU(2,1)$ by angle $\eta$ about the vertical axis
$\{0\}\times\br\subset \sch$.

We can naturally
extend the foliation of the punctured Heisenberg
group $\sch\bs\{0\}$ by Heisenberg spheres $S(0,r)$ to a foliation of the hyperbolic 
space $\ch 2$ by bisectors $\gtS_r$ having those $S(0,r)$ as the spheres at
infinity. It is well known (see \cite{M2}) that each bisector $\gtS_r$
contains a geodesic $\ga_r$ which connects points $(0,-r^2)$ and
$(0,r^2)$ of the Heisenberg group $\sch$ at infinity, and furthermore 
$\gtS_r$ fibers over  $\ga_r$  by complex geodesics $Y$ whose circles at
infinity are complex circles foliating the sphere $S(0,r)$. 

Using those foliations of the hyperbolic 
space $\ch 2$ and bisectors $\gtS_r$, we extend the elementary bending homeomorphism 
$\vp_{\eta,\zeta}\col\ov{\sch}\ra\ov{\sch}$ at infinity to an elementary
bending homeomorphism $\vP_{\eta,\zeta}\col\ov{\ch 2}\ra\ov{\ch 2}$.
Namely, the map $\vP_{\eta,\zeta}$ preserves each of bisectors $\gtS_r$,
each complex geodesic fiber $Y$ in such bisectors, and fixes the intersection points
$y$ of those complex geodesic fibers and the complex geodesic connecting the origin and 
$\infty$ of the Heisenberg group $\sch$ at infinity. We complete our extension
$\vP_{\eta,\zeta}$ by defining its restriction to a given (invariant) complex
geodesic fiber $Y$ with the fixed point $y\in Y$. This map is obtained by  
radiating the circle homeomorphism $\vp_{\eta,\zeta}|_{\p Y}$
to the whole (Poincar\'e) hyperbolic 2-plane $Y$ along geodesic rays
$[y,\infty)\subset Y$, so that it preserves circles in $Y$ centered at $y$ and 
bends (at $y$, by the angle $\eta$) the geodesic in $Y$ connecting the central
points of the corresponding arcs of the complex circle $\p Y$, see
Fig.6.

\midinsert
$$\hss\psfig{file=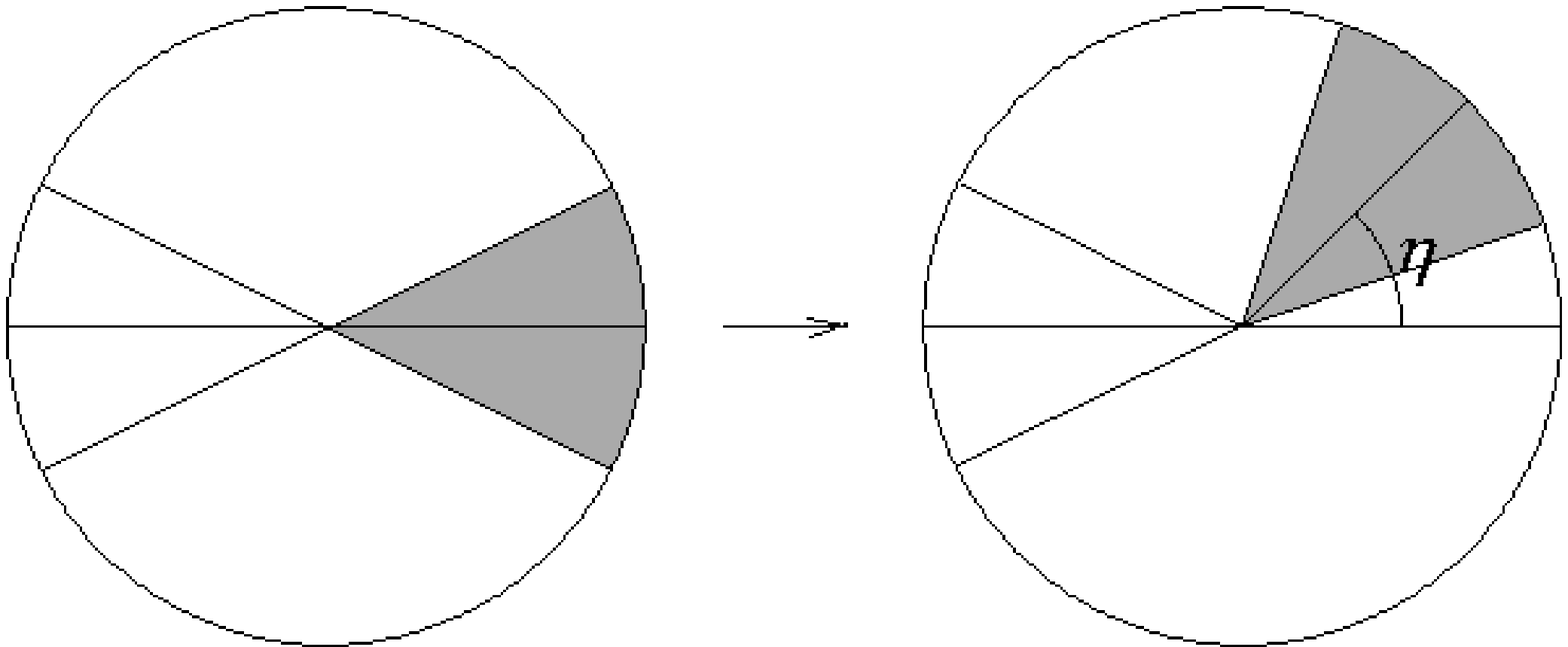,width=400pt}\hss$$
\captionwidth{7pc}
\botcaption{Figure 6} % intersections in figure 5 with round disks
\endcaption
\endinsert

Due to the construction, the elementary bending (quasiconformal) 
homeomorphism $\vP_{\eta,\zeta}$ commutes with elements of the cyclic loxodromic
group $G_0\subset G$. Another most important property of the 
homeomorphism $\vP_{\eta,\zeta}$ is the following. 

Let $D_z(G)$ be 
the Dirichlet fundamental polyhedron of the group $G$ centered at a given
point $z$ on the axis $A$ of the cyclic loxodromic group $G_0\subset G$,
and $\gtS^+\subset \ch 2$ be a ``half-space" disjoint from  $D_z(G)$ and 
bounded by a bisector $\gtS\subset \ch 2$ which is different from bisectors 
$\gtS_r, r>0$, and contains a side $\gts$ of the polyhedron $D_z(G)$. Then there is
an open neighborhood $U(\ov{\gtS^+})\subset \ov{\ch 2}$ such that the restriction
of the elementary bending homeomorphism $\vP_{\eta,\zeta}$ to it either
is the identity or coincides with the unitary rotation $U_{\eta}\subset PU(2,1)$
by the angle $\eta$ about the ``vertical" complex geodesic (containing 
the vertical axis $\{0\}\times\br\subset \sch$ at infinity).

The above properties of quasiconformal homeomorphism $\vP=\vP_{\eta,\zeta}$
show that the image $D_{\eta}=\vP_{\eta,\zeta}(D_z(G))$ is a polyhedron in 
$\ch 2$ bounded by bisectors. Furthermore, there is a natural identification of
its sides induced by $\vP_{\eta,\zeta}$. Namely, the pairs of sides preserved by
$\vP$ are identified by the original generators of the group $G_1\subset G$. 
For other sides $\gts_{\eta}$ of $D_{\eta}$, which are images of corresponding 
sides $\gts\subset D_z(G)$
under the unitary rotation $U_{\eta}$, we define side pairings by using the group
$G$ decomposition (see Fig. 1-4). 

Actually, if $G=G_1 *_{G_0}G_2$,
we change the original side pairings $g\in G_2$ of $D_z(G)$-sides to the 
hyperbolic isometries $U_{\eta}gU_{\eta}^{-1}\in PU(2,1)$. In the case
of HNN-extension, $G=G_1 *_{G_0}=\langle G_1, g_2\rangle$, 
we change the original side pairing $g_2\in G$ of $D_z(G)$-sides to the
hyperbolic isometry $U_{\eta}g_2\in PU(2,1)$. In other words, we define
deformed groups $G_{\eta}\subset PU(2,1)$  correspondingly as

$$G_{\eta}=G_1 *_{G_0}U_{\eta}G_2U^{-1}_{\eta}\quad \text{or}\quad
    G_{\eta}=\langle G_1, U_{\eta}g_2\rangle=G_1 *_{G_0}\,.\tag4.5$$

We see from (4.5) that the family of representations $G\ra G_{\eta}\subset PU(2,1)$
does not depend on 
angles $\zeta$ and holomorphically depends on the angle parameter $\eta$. 
Let us also observe that, for small enough angles $\eta$, the behavior of neighboring 
polyhedra $g'(D_{\eta})$,
$g'\in G_{\eta}$ is the same as of those $g(D_z(G))$, $g\in G$, around
the Dirichlet fundamental polyhedron $D_z(G)$. This is because the new polyhedron 
$D_{\eta}\subset \ch 2$ has isometrically the same (tesselations of) 
neighborhoods of its side-intersections 
as $D_z(G)$ had. This implies that the polyhedra $g'(D_{\eta})$,
$g'\in G_{\eta}$, form a tesselation of $\ch 2$ (with non-overlapping interiors). 
Hence the deformed group $G_{\eta}\subset PU(2,1)$ is a discrete group,
and $D_{\eta}$ is its fundamental polyhedron bounded by bisectors.

Using $G$-compatibility of the restriction of the
elementary bending homeomorphism $\vP=\vP_{\eta,\zeta}$ to the closure
$\ov{D_z(G)}\subset \ov{\ch 2}$, we equivariantly extend it from the polyhedron 
$\ov{D_z(G)}$ to the whole space $\ch 2\cup\Om(G)$ accordingly to the $G$-action.
In fact, in terms of the natural isomorphism $\chi\col G\ra G_{\eta}$ which is 
identical on the subgroup $G_1\subset G$, we can write the obtained 
$G$-equivariant homeomorphism 
$F=F_{\eta}\,\col\, \ov{\ch 2}\bs\La(G)\ra\ov{\ch 2}\bs\La(G_{\eta})$ 
in the following form:

$$\aligned F_{\eta}(x)&=\vP_{\eta}(x)\quad\text{for }\quad x\in
\ov{D_z(G)},\\
     F_{\eta}\circ g(x)&=g_{\eta}\circ F_{\eta}(x)
\quad\text{for }\quad 
x\in \ov{\ch 2}\bs\La(G),\,\, g\in G,\,\, g_{\eta}=\chi(g)\in G_{\eta}\,.
\endaligned\tag4.6$$

Due to quasiconformality of $\vP_{\eta}$, the extended $G$-equivariant 
homeomorphism $F_{\eta}$
is quasiconformal. Furthermore, its extension by continuity to the limit
(real) circle $\La(G)$ coincides with the canonical 
equivariant homeomorphism
$f_{\chi}\col \La(G)\ra\La(G_{\eta})$ given by the isomorphism 
Theorem 2.1. Hence we have a $G$-equivariant
quasiconformal self-homeomorphism of the whole space $\ov{\ch 2}$, which we
denote as before by $F_{\eta}$. 

We claim that the family $\{F^*_{\eta}\}$ of representations
$F^*_{\eta}\,\col\, G\to G_{\eta}=F_{\eta}G_2F^{-1}_{\eta}$, 
$\eta\in (-\eta_0,\,\eta_0)$,
defines a nontrivial curve $\Cal B\col (-\eta_0,\,\eta_0)\ra \Cal R_0(G)$
in the variety $\Cal R_0(G)$ of faithful discrete representations of $G$ 
into $PU(2,1)$, which covers a nontrivial curve
in the Teichm\"uller space $\Cal T(G)$ represented by conjugacy classes
$[\Cal B(\eta)]=[F^*_{\eta}]$. We call the constructed  deformation $\Cal B$ 
the bending deformation
of a given lattice $G\subset PO(2,1)\subset PU(2,1)$ along a bending geodesic 
$A\subset \ch 2$ with loxodromic stabilizer $G_0\subset G$. In terms of
manifolds, $\Cal B$ is the bending deformation of a given complex surface
$M=\ch 2/G$ homotopy equivalent to its totally real geodesic surface 
$S_g\subset M$, along a given simple geodesic $\a$. 
Summarizing properties of these bending deformations, we have:

\proclaim{Theorem 4.7} Let $G\subset PO(2,1)\subset PU(2,1)$ be a given
lattice uniformizing a Riemann 2-surface $S_p$ of genus 
$p\geq 2$. Then, for any simple closed geodesic \break
$\a\subset S_p=H^2_{\br}/G$ and a sufficiently small $\eta_0>0$, 
there is a bending deformation \break
$\Cal B_\a\col (-\eta_0,\,\eta_0)\ra \Cal R_0(G)$ of the group $G$ along $\a$,
$\Cal B_\a(\eta)=\rho_{\eta}=F^*_{\eta}$, 
induced by $G$-equivariant quasiconformal homeomorphisms 
$F_{\eta}: \ov{\ch 2} \ra \ov{\ch 2}$. \qed
\endproclaim

Noticing from (4.5) that the bending deformations along disjoint simple closed
geodesics are independent, we obtain:

\proclaim{Corollary 4.8} Let $S_p=\rh 2/G$ be a closed totally real
geodesic surface of genus $p>1$ in a given complex hyperbolic surface
$M=\ch 2/G$, $G\subset PO(2,1)\subset PU(2,1)$. Then there is an
embedding $\pi\circ\Cal B\col B^{2p-2}\hra \Cal T(M)$ of a real $(2p-2)$-ball
into the Teichm\"uller space of $M$,
defined by bending deformations along disjoint closed geodesics in $M$  and
the projection $\pi\col \Cal D(M)\ra \Cal T(M)=\Cal D(M)/PU(2,1)$.
\endproclaim
  
\demo{Proof} First, we observe that our construction works not only in the case of 
complex surfaces $M$ with totally real geodesic surface but as well as in the case of a
homotopy equivalent surface which has a totally real geodesic piece bounded by closed geodesics. 
In such a case, the holonomy group of that piece is a non-elementary discrete 
subgroup in $PU(2,1)$ preserving a totally real geodesic plane 
$\rh 2\subset\ch 2$. Then we see from (4.5) that bendings along disjoint closed geodesics
are independent.
In addition to the above construction, we have also to show that 
our bending deformation is not trivial, and 
$[\Cal B(\eta)]\neq [\Cal B(\eta')]$ for any $\eta\neq\eta'$. 

The non-triviality of a bending deformation follows directly from (4.5), cf. \cite{A5}.
Namely, the restrictions $\rho_{\eta}|_{G_1}$ of bending representations
to a non-elementary subgroup $G_1\subset G$ (in general, to a ``real" subgroup
$G_r\subset G$ corresponding to a totally real geodesic piece in the homotopy
equivalent surface $S\backsimeq M$) are identical. So if the deformation $\Cal B$ 
were trivial then it would be conjugation of the group $G$ by projective
transformations that commute with the non-trivial real subgroup $G_r\subset G$ and 
pointwise fix the totally real geodesic plane $\rh 2$. This contradicts
to the fact that the limit set of any deformed group $G_{\eta}$,
$\eta\neq 0$, does not belong to the real circle containing the Cantor
limit set $\La(G_r)$.

The injectivity of the map $\Cal B$ can be obtained by using \'Elie Cartan \cite{Ca}
angular invariant  $\ba(x)$ for a triple $x=(x^0, x^1, x^2)$ of points in 
$\p\ch 2$. It satisfies the properties (see \cite{G3}): 
\roster
\item $-\pi/2\leq\ba(x)\leq\pi/2$;
\item $\ba(x)=0$ if and only if  $x^0, x^1$ and $x^2$ lie on an $\br$-circle; 
\item $\ba(x)=\pm \pi/2$ if and only if $x^0, x^1$ and $x^2$ lie on a
chain ($\bc$-circle);
\item For two triples $x$ and $y$, $\ba(x)=\ba(y)$ if and only if
there exists $g\in PU(2,1)$ such that $y=g(x)$; furthermore, such a $g$ is unique
provided that $\ba(x)$ is neither zero nor  $\,\pm \pi/2$.
\item  $\ba(x)= - \ba(y)$ if and only if
there exists an antiholomorphic automorphism $g\in \is \ch 2$ such
that $y=g(x)$.
\endroster

 Namely, let $g_2\in G\bs G_1$ be a generator of
the group $G$ in (4.5) whose fixed point $x^2\in \La(G)$ lies in $\br_+\times \{0\}
\subset \sch$, and $x^2_{\eta}\in \La(G_{\eta})$ the corresponding fixed point
of the element $\chi_{\eta}(g_2)\in G_{\eta}$ under the free-product
isomorphism $\chi_{\eta}\col G\ra G_{\eta}$.
Due to our construction, one can see that the orbit
$\ga(x^2_{\eta})$, $\ga\in G_0$, under the loxodromic (dilation) subgroup $G_0\subset 
G\cap G_{\eta}$ approximates the origin along a ray 
$(0,\infty)$ which has a non-zero
angle $\eta$ with the ray $\br_-\times \{0\}\subset \sch$. The latter ray also contains
an orbit $\ga(x^1)$, $\ga\in G_0$, of a limit point $x^1$ of $G_1$ which approximates 
the origin from the other side. Taking triples $x=(x^1, 0, x^2)$ and 
$x_{\eta}=(x^1, 0, x^2_{\eta})$ of points which lie correspondingly in
the limit sets $\La(G)$ and $\La(G_{\eta})$, we have
that $\ba(x)=0$ and $\ba(x_{\eta})\neq 0, \,\pm\pi/2$. Due to Theorem 2.1,
both limit sets are topological circles which however cannot be equivalent
under a hyperbolic isometry because of different Cartan invariants
(and hence, again, our deformation is not
trivial). 

Similarly, for two different values $\eta$ and $\eta'$, we have triples
$x_{\eta}$ and $x_{\eta'}$ with different (non-trivial) Cartan angular invariants 
$\ba(x_{\eta})\neq \ba(x_{\eta'})$. Hence two topological circles 
$\La(G_{\eta})$ and $\La(G_{\eta'})$ are not hyperbolically isometric. This
completes the proof of injectivity and of the whole Corollary. \qed
\enddemo

One can apply the above proof to a general situation of bending
deformations of a complex hyperbolic surface $M=\ch 2/G$ whose holonomy 
group  $G\subset PU(2,1)$ has a non-elementary subgroup $G_r$ preserving
a totally real geodesic plane $\rh 2$. In other words, such a complex surfaces $M$
has an embedded totally real geodesic surface with geodesic boundary.
 In particular all complex
surfaces constructed in \cite{GKL} with a given Toledo invariant lie in
this class. So we immediately have:

\proclaim{Corollary 4.9} Let $M=\ch 2/G$ be a complex hyperbolic surface
with embedded totally real geodesic surface $S_r\subset M$
with geodesic boundary, and $\Cal B\col (-\eta, \eta)\ra \Cal D(M)$
a bending deformation of $M$ along a simple closed geodesic $\a\subset S_r$. 
Then the map $\pi\circ\Cal B\col (-\eta, \eta)\ra \Cal T(M)=\Cal D(M)/PU(2,1)$ is a smooth
embedding provided that the limit set $\La(G)$ of the holonomy group $G$
does not belong to the $G$-orbit of the real circle $S^1_{\br}$ and the
chain $S^1_{\bc}$, where the latter is the infinity of the complex geodesic containing 
a lift $\tilde\a\subset\ch 2$ of the closed geodesic $\a$, and the former one contains
the limit set of the holonomy group $G_r\subset G$ of the geodesic
surface $S_r$. \qed
\endproclaim

\remark{Remark 4.10} It follows from the above construction of the 
bending homeomorphism $F_{\eta, \zeta}$, that the deformed complex
hyperbolic surface $M_{\eta}=\ch 2/G_{\eta}$ fibers over the pleated 
hyperbolic surface $S_{\eta}=F_{\eta}(\rh 2)/G_{\eta}$ (with the closed geodesic 
$\a$ as the singular locus). The fibers of this fibration are
``singular real planes" obtained from totally real geodesic 2-planes by
bending them  by angle $\eta$ along complete real geodesics. These
(singular) real geodesics are the intersections of the complex geodesic
connecting the axis $A$ of the cyclic group $G_0\subset G$ and the
totally real geodesic planes that represent fibers of the original fibration in
$M=\ch 2/G$.
\endremark
\vfil
\eject

\def\ref#1{[#1]}
\eightpoint
\parindent=36pt

\head  REFERENCES\endhead
\bigskip

\frenchspacing

\item{\ref{A1}}  Boris Apanasov, Discrete groups in Space and
Uniformization  Problems. - Math. and Appl., {\bf 40}, Kluwer
Academic Publishers, Dordrecht, 1991.

\item{\ref{A2}}  \underbar{\phantom{Apanasov}}, Nontriviality  of Teichm\"uller space for
Kleinian  group in  space.-  Riemann Surfaces and Related
Topics: Proc. 1978 Stony  Brook Conference (I.Kra and
B.Maskit, eds), Ann. of  Math. Studies {\bf 97},
Princeton Univ. Press, 1981, 21-31.

\item{\ref{A3}}    \underbar{\phantom{Apanasov}},   Thurston's bends and geometric
deformations of conformal structures.- Complex Analysis and
Applications'85,  Publ. Bulgarian Acad. Sci.,
Sofia, 1986, 14-28.

\item{\ref{A4}}  \underbar{\phantom{Apanasov}}, Deformations of conformal structures on
hyperbolic manifolds.- J. Diff. Geom. {\bf 35}
(1992), 1-20.

\item{\ref{A5}} \underbar{\phantom{Apanasov}}, Conformal geometry of
discrete groups and manifolds.  - W. de Gruyter, Berlin- New York, to appear.

\item{\ref{A6}} \underbar{\phantom{Apanasov}}, Canonical homeomorphisms in
Heisenberg groups induced by isomorphisms of discrete subgroups of $PU(n,1)$.
 - Russian Acad. Sci. Dokl. Math., to appear.

\item{\ref{A7}} \underbar{\phantom{Apanasov}}, Quasiconformality and
geometrical finiteness in Carnot-Carath\'eodory and negatively curved
spaces. - Math. Sci. Res. Inst. at Berkeley, 1996-019.

\item{\ref{ACG}} Boris Apanasov, Mario Carneiro and Nikolay Gusevskii,
Some deformations of complex hyperbolic surfaces. - In preparation.
 
\item{\ref{AG}} Boris Apanasov and Nikolay Gusevskii, 
The boundary of Teichm\"uller space of complex hyperbolic surfaces. - In preparation.

\item{\ref{AX1}}  Boris Apanasov and Xiangdong Xie, Geometrically finite 
complex hyperbolic manifolds.- Preprint, 1995.

\item{\ref{AX2}} \underbar{\phantom{Apanasov}}, Manifolds of negative
curvature and nilpotent groups.- Preprint, 1995.

\item{\ref{Be}} A.F.~Beardon, The geometry of discrete groups.- Springer-Verlag,
Berlin/New York, 1983.

\item{\ref{Bu}} Peter Buser, The collar theorem and examples.- Manuscr.
Math. {\bf 25} (1978), 349-357.

%\item{\ref{Be}} Igor Belegradek, Discrete surface groups actions with
%accidental parabolics on complex hyperbolic plane.- Preprint, 1995.

\item{\ref{BS}} D.~Burns and S.~Shnider, Spherical hypersurfaces in complex
manifolds. - Inv. Math., {\bf 33}(1976), 223-246.

\item{\ref{CG}} S.~Chen and L.~Greenberg, Hyperbolic spaces.-
Contributions to Analysis, Academic Press, New York, 1974, 49-87.

\item{\ref{Ca}} \'Elie Cartan, Sur le groupe de la g\'eom\'etrie
hypersph\'erique.- Comm. Math. Helv. {\bf 4} (1932), 158-171.

\item{\ref{C1}} Kevin Corlette, %Hausdorff dimensions of limit sets. -
%Invent. Math. {\bf 102} (1990), 521-542.
%\item{\ref{C2}}  \underbar{\phantom{Apanasov}}, 
Archimedian superrigidity and hyperbolic geometry. - 
Ann. of Math. {\bf 135} (1992), 165-182.

\item{\ref{Cy}} J.~Cygan, Wiener's test for Brownian motion on the Heisenberg
group. - Colloquium Math. {\bf 39}(1978), 367-373.

\item{\ref{FG}} Elisha Falbel and Nikolay Gusevskii, Spherical CR-manifolds of
dimension 3.- Bol. Soc.Bras.Mat. {\bf 25}(1994), 31-56

\item{\ref{G1}} William Goldman, Representations of fundamental groups of
surfaces.- Geometry and topology (J.Alexander and J.Harer, Eds), Lect.
Notes Math. {\bf 1167}, Springer, 1985, 95-117. 

\item{\ref{G2}} \underbar{\phantom{Apanasov}}, Geometric structures on
manifolds and varieties of representations. - Geometry of Group Representations
(A.Magid and W. Goldman, Eds), Contemp. Math. {\bf 74} (1988), 169-198.
 
\item{\ref{G3}} \underbar{\phantom{Apanasov}}, Complex hyperbolic
geometry. - Oxford Univ. Press, to appear.

\item{\ref{GKL}} W.Goldman, M.Kapovich and B.Leeb, Complex hyperbolic
manifolds homotopy equivalent to a Riemann surface. - Preprint, 1995.

\item{\ref{GM}} William Goldman and John Millson, Local rigidity of discrete
groups acting on complex hyperbolic space.- Invent. Math. {\bf 88}(1987),
495-520.

\item{\ref{GP}} William Goldman and John Parker,  Complex hyperbolic ideal
triangle groups.- J. reine angew. Math. {\bf 425} (1992), 71-86.

\item{\ref{Gu}} Nikolay Gusevskii, Colloquium talk, Univ. of Oklahoma,
Norman, December 1995.

\item{\ref{JM}}  Dennis Johnson and  John Millson,  Deformation  spaces
associated  to  compact  hyperbolic  manifolds.- Discrete
Groups in Geometry  and Analysis:  Papers in  Honor of G.D.
Mostow on  His Sixtieth Birthday,  Ed. R. Howe, Birkhauser,
Boston, 1987, 48-106.

\item{\ref{Ke}} Linda Keen, Collars on Riemann surfaces.- Discontinuous
groups and Riemann surfaces, Ann. of Math. St. {\bf 79}, Princeton Univ.
Press, 1974, 263-268.

\item{\ref{KR}} Adam Koranyi and Martin Reimann, Quasiconformal mappings on the
Heisenberg group. - Invent. Math. {\bf 80} (1985), 309-338.

\item{\ref{Ko}} Christos Kourouniotis,
Deformations of hyperbolic structures.- Math. Proc. Cambr. Phil. Soc.
 {\bf 98} (1985), 247-261.

\item{\ref{Mas}} Bernard Maskit, Kleinian groups. - Springer-Verlag, 1987.

\item{\ref{Ma}} Georgii A. Margulis, Discrete groups of motions of
manifolds of non-positive curvature.- Proc. ICM, Vancouver, v.2, 1974, 21-34 (Russian);
Engl. Transl.: Amer. Math. Soc. Translations {\bf 109} (1977), 33-45.

\item{\ref{M1}} George D. Mostow, Strong rigidity of locally symmetric
spaces.- Princeton Univ. Press, 1973.

\item{\ref{M2}} \underbar{\phantom{Apanasov}},  On a remarkable class of polyhedra in
complex hyperbolic space. - Pacific J. Math., {\bf 86}(1980), 171-276.

\item{\ref{P}} Pierre Pansu, M\'etriques de Carnot-Carath\'eodory et
quasiisom\'etries des espaces symm\'etries de rang un.- Ann. Math. 
{\bf 129} (1989), 1-60.

\item{\ref{Pr}} John Parker, Shimizu's lemma for complex hyperbolic space. -
Intern. J. Math. {\bf 3:2} (1992), 291-308.

\item{\ref{Th}} William Thurston, The geometry and topology of
three-manifolds. - Lect. Notes, Princeton Univ., 1981.

\item{\ref{To}} Domingo Toledo, Representations of surface groups on complex
hyperbolic space.- J. Diff. Geom. {\bf 29} (1989), 125-133.

\item{\ref{Tu}} Pekka Tukia, On isomorphisms of geometrically finite
Kleinian groups.- Publ. Math. IHES {\bf 61}(1985), 171-214.

\item{\ref{V}} Serguei K. Vodopyanov, Quasiconformal mappings on Carnot
groups.- Russian Dokl. Math., to appear

\item{\ref{Ya}} S.T.Yau, Calabi's conjecture and some new results in
algebraic geometry.- Proc. Nat. Acad. Sci {\bf 74} (1977), 1798-1799.

\item{\ref{Yu}} Chengbo Yue, Dimension and rigidity of quasi-Fuchsian
representations.- Ann. of Math. {\bf 143} (1996)331-355.

\enddocument